\documentclass[11pt]{article}

\usepackage{arxiv_math,latexsym}

\newcommand{\B}{\mathcal{B}}

\def\elead{\mbox{\sf ELead}}
\def\lead{\mbox{\sf Lead}}
\def\mon{\mbox{\sf Mon}}

\begin{document}

\title{Balls-in-bins with feedback and Brownian Motion}
\author{Roberto Oliveira\thanks{IBM T.J. Watson Research Center, Yorktown Heights, NY 10598. \texttt{riolivei@us.ibm.com}. Work mostly done while the author was a Ph.D. student at New York University under a CNPq scholarship.}} \maketitle
\begin{abstract}In a balls-in-bins process with feedback, balls are
sequentially thrown into bins so that the probability that a bin
with $n$ balls obtains the next ball is proportional to $f(n)$ for
some function $f$. A commonly studied case where there are two
bins and $f(n) = n^p$ for $p > 0$, and our goal is to study the
fine behavior of this process with two bins and a large initial
number $t$ of balls. Perhaps surprisingly, Brownian Motions are an essential part of both our proofs.\\
For $p>1/2$, it was known that with probability $1$ one of the
bins will lead the process at all large enough times. We show that
if the first bin starts with $t+\lambda\sqrt{t}$ balls (for
constant $\lambda\in \R$), the probability that it always or
eventually leads has a non-trivial limit depending on $\lambda$.\\
For $p\leq 1/2$, it was known that with probability $1$ the bins
will alternate in leadership. We show, however, that if the
initial fraction of balls in one of the bins is $>1/2$, the time
until it is overtaken by the remaining bin scales like
$\Theta({t^{1+1/(1-2p)}})$ for $p<1/2$ and $\exp(\theita{t})$ for
$p=1/2$. In fact, the overtaking time has a non-trivial distribution around the scaling factors, which we determine explicitly.\\
Our proofs use a continuous-time embedding of the balls-in-bins
process (due to Rubin) and a non-standard approximation of the
process by Brownian Motion. The techniques presented also extend
to more general functions $f$.\end{abstract}

\section{Introduction}\label{sec:intro}

\subsection{The process}

If $f:\N\to (0,+\infty)$ is a positive function, a {\em
balls-in-bins process with feedback function $f$} a discrete-time
Markov process with $B$ bins, each one of which containing
$I_i(m)>0$ balls at time $m$ for each $m\in \{0,1,2,\dots\}$ and
$i\in\{1,\dots,B\}$. Its evolution is as follows: at each time
$m>0$, a ball is added to a bin $i_m$, so that
$I_{i_m}(m)=I_{i_m}(m-1)+1$ and $I_j(m)=I_j(m-1)$ for all
$i\in\{1,\dots,B\}\backslash\{i_m\}$, and the random choice of bin
$i_m$ has distribution
\begin{equation}\label{eq:recipe}\Pr{i_{m}=i\mid \{I_j(m-1)\,:\,{1\leq j\leq B}\}}=
\frac{f(I_i(m-1))}{\sum_{j=1}^B f(I_j(m-1))}\;\; (1\leq i\leq
B),\end{equation}We will commonly refer to this recipe by saying
that {\em bin $i$ receives a ball at time $m$ (i.e. $i_m=i$) with
probability proportional to $f(I_m(i))$}.

Such processes\footnote{A longer background discussion is
available from the author's PhD thesis \cite{tese}.} were
introduced to the Discrete Mathematics community by Drinea, Frieze
and Mitzenmacher \cite{Drinea02}, where they were motivated by
economical problems of competition and mathematically related
preferential attachment models for large networks\footnote{More
specifically, the model introduced by Krapivsky and Redner
\cite{Krapivsky02} and independently by Driena, Enachescu and
Mitzemacher \cite{DrineaEM}. This model generalizes the
Barab\'{a}si-\'{A}lbert model, which is discussed in the survey
\cite{AlbertSurvey}.}. That paper treats only the case where
$f(x)=x^p$ for some parameter $p>0$. In this case $f$ is
increasing, and therefore {\em the rich get richer}: the more
balls a bin has, the more likely it is to receive the next ball.
In economic terms, one could think of bins as competing products
and balls as customers; in that case, the more popular a product
is, the more likely it is to obtain a new customer.

The main question addressed in that paper is whether this
phenomenon ensures large advantages {\em in the long run} for some
bin. The authors show that if $p>1$, there almost surely exists
one bin that gets all but a negligible fraction of the balls in
the large-time limit; whereas for $p<1$, the asymptotic fractions
of balls in each bin even out with time. The $p=1$ case is the
classic P\'olya Urn model, for which it has been long known that
the number of balls in each bin converges almost surely to a
non-degenerate random variable, and thus the process has different
{\em regimes} depending on the choice of parameter $p$.

\subsection{The three regimes}

However, a much stronger result is available. A paper by Khanin
and Khanin \cite{Khanin01} introduced what amounts to the same
process as a model for neuron growth, and proved that if $p>1$,
there almost surely is some bin that gets {\em all but finitely
many balls}, an event that we call {\em monopoly}. They also show
that for $1/2<p\leq 1$, monopoly has probability $0$, but there
almost surely will be some bin which will lead the process from
some finite time on (we call this {\em eventual leadership}),
whereas this cannot happen if $0<p\leq 1/2$. Therefore, the
balls-in-bins process has {\em three regimes} of behavior
corresponding to three ranges of $p$.

In fact, the result of \cite{Khanin01} generalizes to any $f$ with
$\min_{x\in\N}f(x)>0$. Consider the following events, which we
call {\em monopoly by bin $i\in[B]$},
\begin{eqnarray}\label{eq:defmon}\mon_i &\equiv& \{\exists M\in \N\;\forall m\geq M\;\;\forall j\in [B]\;\;j\neq i\Rightarrow I_m(j)=I_M(j)\}\\ \nonumber &=&\{\exists M\in \N\;\forall m\geq M\;\;i_{m}=i\};\end{eqnarray}
and {\em eventual leadership by bin $i\in[B]$}:
\begin{equation}\label{eq:defelead}\elead_i \equiv \{\exists M\in \N\;\forall m\geq M\;\;\forall j\in [B]\;\;j\neq i\Rightarrow I_m(j)<I_m(i)\}.\end{equation}
Clearly, $\elead_i\supset \mon_i$. A not-too-difficult extension
of \cite{Khanin01} (proven in \cite{Spencer??,tese}) says that
\begin{theorem}[From \cite{Khanin01,Spencer??,tese}]\label{thm:Khanin}If
$\{I_m\}_{m=0}^{+\infty}$ is a balls-in-bins process with $B$ bins
and feedback function $f=f(x)\geq c$ for some $c>0$, then there
are three mutually exclusive possibilities.
\begin{enumerate}
\item $\sum_{n\geq 1}f(n)^{-1} < +\infty$, in which case
$\Pr{\cup_{i=1}^B\,\mon_i} = \Pr{\cup_{i=1}^B\,\elead_i}=1$ (we
call this the {\em monopolistic} regime); \item $\sum_{n\geq
1}f(n)^{-1} = +\infty$ but $\sum_{n\geq 1}f(n)^{-2} < +\infty$, in
which case $\Pr{\cup_{i=1}^B\,\mon_i} = 0$ but
$\Pr{\cup_{i=1}^B\,\elead_i}=1$ (this is the {\em eventual
leadership regime}); \item $\sum_{n\geq 1}f(n)^{-2} =+\infty$, in
which case $\Pr{\cup_{i=1}^B\,\mon_i} =
\Pr{\cup_{i=1}^B\,\elead_i}=0$ (this is the {\em almost-balanced
regime}).
\end{enumerate}
This holds irrespective of the initial conditions of the
process.\end{theorem} Notice that the three cases of the Theorem
applied to the $f(x)=x^p$ family correspond to $p>1$, $1/2<p\leq
1$ and $0\leq p<1/2$; in other words, this family of $f$ has {\em
phase transitions} at $p=1$ and $p=1/2$. We sketch a proof of this
result in \secref{exp_proof}, both for completeness and to give
readers a better acquaintance with the techniques in the present
paper.

\subsection{The present work}

This paper is part of a series by the present author in
collaboration with Michael Mitzenmacher and Joel Spencer
\cite{MitzenmacherOS04,Oliveira05,avoid} that attempts to quantify
different aspects of the three qualitative regimes presented in
\thmref{Khanin}. Our specific purpose in the present paper is to
prove two not-quite-related results about these processes in
different regimes, when the initial number of balls on both bins
is large. What brings these two results is that both proofs use
Brownian Motion in an unexpected and surprising way.

Our first result is a {\em scaling result for the eventual
leadership and monopoly regimes}. Suppose, for simplicity, that
$f(x)=x^p$ with $p>1/2$. Recall the definition of eventual
leadership by bin $i$ \eqnref{defelead}, and let $\lead_i$ be the
event that bin $1$ leads the process at {\em all times}:
\begin{equation}\label{eq:deflead}\lead_i \equiv \{\forall m\geq 0\;\; \forall j\in[B]\;\; j\neq i\Rightarrow I_m(j)<I_m(i)\}.\end{equation}

If the initial number of bins is $t\gg 1$ and $I_0(1)\approx t/2$,
then $\Pr{\lead_i}\approx 0$ and $\Pr{\elead_i}\approx 1$ for
$i=1,2$. On the other hand, if $I_0(1)$ is much larger than
$I_0(2)$, one think that $\Pr{\lead_1},\Pr{\elead_1}\approx 1$.
Our question is {\em how large is large enough}? That is to say,
at what {\em scale} do these two probabilities grow from $0$ to
$1$? We show that the answer is in fact $\theita{\sqrt{t}}$, and
give an exact asymptotic result.

\begin{theorem}\label{thm:scaling_cor}Let $\{\lambda_t\}_{t\in\N}\subset\R$ form a sequence such that such that
$$\lambda\equiv \lim_{t\to +\infty}\lambda_t \in\R\mbox{ exists}$$
and
$$\forall t\in\N,\;\;\frac{t}{2}\pm \lambda_t\sqrt{\frac{t}{4p-2}}\in \N.$$
Assume that for each $t$, $\{I^{(t)}_m(j)\,:\,j\in[2],m\geq 0\}$
is a balls-in-bins process with two bins, feedback function
$f(x)=x^p$, $1/2<p<+\infty$ and initial conditions
$$(x(t),y(t))\equiv
\left(\frac{t}{2}+\lambda_t\sqrt{\frac{t}{4p-2}},\frac{t}{2}-\lambda_t\sqrt{\frac{t}{4p-2}}\right).$$
Then
\begin{eqnarray}\lim_{t\to +\infty}\Prp{(x(t),y(t))}{\elead_1} &=& \gamma(\lambda)\equiv \sqrt{\frac{1}{2\pi}}\int_{-\infty}^{\lambda}e^{-x^2/2}\,dx,\\
\lim_{t\to +\infty}\Prp{(x(t),y(t))}{\lead_1} &=&
\Gamma(\lambda)\equiv\sqrt{\frac{2}{\pi}}\int_{0}^{\max\{\lambda,0\}}e^{-x^2/2}\,dx.\end{eqnarray}
\end{theorem}

This theorem is an extension of a result by Mitzenmacher, Oliveira
and Spencer \cite{MitzenmacherOS04}, who showed a similar scaling
for $\mon_1$ when $p>1$. That paper used Ess\'{e}en's inequality
for approximation by Gaussians together with a continuous-time
embedding of the balls-in-bins process; we shall also use the
latter device together with approximation by Brownian Motion,
especially to estimate $\Prp{(x(t),y(t))}{\mon_i}$. Notice that as
$p\searrow 1/2$ the scaling term $\sqrt{\frac{t}{4p-2}}$ becomes
bigger; i.e. near the $p=1/2$ phase transition, it becomes harder
to bias the process towards (eventual) leadership by either bin.

The second result we prove is about the almost-balanced case.
Suppose, again for simplicity, that there are two bins ($B=2$) and
$f(x)=x^p$, $0<p\leq 1/2$, . In this case \thmref{Khanin} says
that for any initial conditions $I_0(1),I_0(2)$ with
$I_0(1)<I_0(2)$, there is a time $m\geq 0$ such that
$I_m(1)>I_m(2)$. Call the first such time the {\em overtaking
time} $V$. By the above, $V<+\infty$, but we have no idea of the
distribution of $V$, and thus we don't know how long the
overtaking might take. We show that if the initial number of balls
is large and bin $2$ has a non-negligibly bigger fraction of the
initial balls, then $V$ can actually be quite large; moreover, it
has an explicit asymptotic distribution.

\begin{theorem}\label{thm:almost_cor}Let $V_{t,\alpha}$ be the overtaking time in a balls-in-bins process with feedback function $f(x)=x^p$ (with $p\in(0,1/2]$ constant) and initial conditions $(\lceil \alpha t\rceil,t-\lceil \alpha t\rceil)$ for $0<\alpha<1/2$. Then there exist random variables $\{U_{t,\alpha,p}\}_{t\in\N}$ such
that
\begin{equation}V_{t,\alpha} =\left\{\begin{array}{l} 2\left\{\frac{(1-2p)}{(1-p)^2}[(1-\alpha)^{1-p} -\alpha^{1-p}]^2 + \frac{(1-\alpha)^{1-2p}\,U_{t,\alpha,p}^2}{t}\right\}^{\frac{1}{1-2p}}\frac{t^{1+\frac{1}{1-2p}}}{U_{t,\alpha,p}^{1+\frac{2p}{1-2p}}} - (t+1)\\
\;\;\;\;\;\;\;\;\;\;\;\;\;\mbox{ if }0<p<\frac{1}{2}\commaeq\\
2(1-\alpha)\, t\,
\exp\left\{4[1-2\sqrt{\alpha(1-\alpha)}]\,\frac{t}{U_{t,\alpha,\frac{1}{2}}^2}\right\}
-(t+1)\\ \;\;\;\;\;\;\;\;\;\;\;\;\;\mbox{ if }p
=\frac{1}{2}\commaeq\end{array}\right.\end{equation} with
probability tending to $1$ as $t\to +\infty$, and
$$U_{t,\alpha,p} \weakto |N|\commaeq$$where $N$ is a standard Gaussian
random variable.\end{theorem}

This means that $V$ becomes larger and larger as $p\nearrow 1/2$,
culminating with the exponential behavior at the phase transition
point $p=1/2$. The economically-inclined might wish to deduce from
this theorem that, under appropriate initial conditions, a
product's leadership might last a long time even in markets with
no propensity for breeding monopolies or ``eternal leaders".

\subsection{Techniques and outline}

Our results in this paper are actually more general: they extend
to a broader (though not entirely general) class of functions $f$
and, in the case of \thmref{scaling_main}, to more than two bins.
All proofs below are done for this more general case and then
specialized for $f(x)=x^p$.

Our proofs have their first two steps in common. The first step
has been employed in \cite{Khanin01,Spencer??} and other works,
and seems to have originated in Davis' work on reinforced random
walks \cite{Davis90}. We shall {\em embed} the discrete-time
process we are interested in into a continuous-time process built
from exponentially distributed random variables, so that
interarrival times at different bins are independent and have an
explicit distribution, which is very helpful in calculations. We
call this the {\em exponential embedding} of the process.

The second technique we use is {\em approximation by Brownian
motion} via Donsker's Invariance Principle. While neither
technique is novel, their conjunction in the way presented here
yields surprising explicit results in the asymptotic regime, once
the appropriate calculations are done.

The remainder of the paper is organized as follows. We discuss
preliminary material in \secref{prelim}. \secref{exp_embed}
rigorously introduces the exponential embedding process and
discusses its key properties. In \secref{technical} we detail the
assumptions we make on our feedback functions $f$, while also
deriving some consequences of those assumptions. \secref{scaling}
proves the general version of \thmref{scaling_cor}, whereas
\secref{almost} contains the proof of the generalization of
\thmref{almost_cor}. Finally, \secref{last} discusses extensions
to our results and some open problems.

\subsection{Acknowledgements}

Michael Mitzenmacher and my Ph.D. advisor Joel Spencer introduced
me to this topic and stimulated me with several interesting
questions and feedback. I also thank Eleni Drinea for useful
discussions.

\section{Preliminaries}\label{sec:prelim}

\par {\em General notation.} Throughout the paper,
$\N=\{1,2,3,\dots\}$ is the set of non-negative integers, $\R^+ =
[0,+\infty)$ is the set of non-negative reals, and for any
$k\in\N\backslash\{0\}$ $[k]=\{1,\dots,k\}$.\\

\par {\em Asymptotics.} We use the standard $O/o/\Omega/\Theta$ notation. The expressions ``$a_n\sim b_n$ as $n\to n_0$" and ``$a_n\ll b_n$ as $n\to n_0$" mean that $\lim_{n\to n_0}(a_n/b_n) = 1$ and $\lim_{n\to n_0}(a_n/b_n) = 0$, respectively.\\

\par {\em Balls-in-bins.} Formally, a balls-in-bins process with feedback function $f:\N\to (0+\infty)$ and $B\in\N$ bins is a discrete-time Markov chain $\{(I_1(m),\dots,I_B(m))\}_{m=0}^{+\infty}$ with state space $\N^B$ and transition probabilities given in the Introduction (see \eqnref{recipe}). We will usually refer to the index $i_m\in[B]$ as {\em the bin that receives a ball at time $m$}. For any $B$, if $E$ is an event of the process and $u\in\N^B$,
$\Prp{u}{E}$ is the probability of $E$ when the initial conditions
are set to $u$.\\

\par {\em Exponential random variables.} $X\eqdist \exp(\lambda)$
means that $X$ is a random variable with exponential distribution
with rate $\lambda>0$, meaning that $X\geq 0$ and
$$\Pr{X>t} = e^{-\lambda t}\;\; (t\geq 0).$$
The shorthand $\exp(\lambda)$ will also denote a generic random
variable with that distribution. Some elementary but extremely
useful properties of those random variables include
\begin{enumerate}\item {\em Lack of memory.} Let $X\eqdist \exp(\lambda)$
and $Z\geq 0$ be independent from $X$. The distribution of $X-Z$
conditioned on $X>Z$ is still equal to $\exp(\lambda)$. \item {\em
Minimum property.} Let $\{X_i\eqdist \exp(\lambda_i)\}_{i=1}^{m}$
be independent. Then $$X_{\min}\equiv \min_{1\leq i\leq
m}X_i\eqdist \exp(\lambda_1+\lambda_2+\dots \lambda_m)$$ and for
all $1\leq i\leq m$
\begin{equation}\Pr{X_i = X_{\min}} = \frac{\lambda_i}{\lambda_1+\lambda_2+\dots
\lambda_m}\end{equation}  \item {\em Multiplication property.} If
$X\eqdist\exp(\lambda)$ and $\eta>0$ is a fixed number, $\eta X
\eqdist\exp(\lambda/\eta)$. \item {\em Moments and transforms.} If
$X\eqdist\exp(\lambda)$, $r\in\N$ and $t\in\R$,
\begin{eqnarray}
\Ex{X^r} & = & \frac{r!}{\lambda^r}\commaeq \\
\Ex{e^{t X}} & = & \left\{\begin{array}{ll}\frac{1}{1 - \frac{t}{\lambda}} & (t<\lambda)\\
                    +\infty  & (t\geq
                    \lambda)\end{array}\right.\end{eqnarray}

\end{enumerate}

\par {\em Weak convergence.} $X_n\weakto Y$ means that the sequence
$\{X_n\}$ of random variables converges weakly to $Y$ as $n\to
+\infty$.

\par {\em Gaussians and cumulative distribution functions.} Finally, we restate the definitions of $\Phi$ and $\Gamma$ in \thmref{scaling_main}:

\begin{eqnarray}\label{eq:defphi}\Phi(\lambda)&\equiv &\sqrt{\frac{1}{2\pi}}\int_{-\infty}^{\lambda}e^{-x^2/2}\,dx,\\
\label{eq:defgamma}\Gamma(\lambda)&\equiv&\sqrt{\frac{2}{\pi}}\int_{0}^{\max\{\lambda,0\}}e^{-x^2/2}\,dx.\end{eqnarray}

If $N$ is a standard Gaussian random variable, $\Gamma$ is the cdf
(cumulative distribution function) of $N$ and $\Gamma$ is the cdf
of $|N|$.

\section{The exponential embedding}\label{sec:exp_embed}

\subsection{Definition and key properties}

Let $f:\N\to(0,+\infty)$ be a function, $B\in\N$ and
$(a_1,\dots,a_B)\in\N^B$. We define below a continuous-time
process with state space $(\N\cup\{+\infty\})^B$ and initial state
$(a_1,\dots,a_B)$ as follows. Consider a set $\{X(i,j)\,:\,i\in
[B],\, j\in \N\}$ of independent random variables, with
$X(i,j)\eqdist \exp(f(j))$ for all $(i,j)\in [B]\times \N$, and
define

\begin{equation}\label{eq:exp_contproc}N_i(t)\equiv \sup\left\{n\in
\N\,:\,\sum_{j=a_i}^{n-1} X(i,j) \leq t\right\}\;\;\;(i\in
[B],t\in\R^+ = [0,+\infty))\commaeq\end{equation} where by
definition $\sum_{j=i}^{k}(\dots)=0$ if $i>k$. Thus $N_i(0)=a_i$
for each $i\in [B]$, and one could well have $N_i(T)=+\infty$ for
some finite time $T$ (indeed, that {\em will} happen for our cases
of interest); but in any case, the above defines a continuous-time
stochastic process, and in fact the $\{N_i(\cdot)\}_{i=1}^B$
processes are independent. Each one of this processes is said to
correspond to {\em bin} $i$, and each one of the times
$$X(i,a_i),X(i,a_{i})+X(i,a_{i}+1),X(i,a_{i})+X(i,a_{i}+1)+X(i,a_{i}+2),\dots$$
is said to be an {\em arrival time at bin $i$}. As in the
balls-in-bins process, we imagine that each arrival correspond to
a ball being placed in bin $i$.

In fact, we {\em claim} that this process is related as follows to
the balls-in-bins process with feedback function $f$, $B$ bins and
initial conditions $(a_1,\dots,a_B)$.

\begin{theorem}[Proven in \cite{Davis90,Khanin01,Spencer??,tese,Oliveira05}]\label{thm:exp_embed}Let the $\{N_i(\cdot)\}_{i\in[B]}$ process be defined as above. One can order the arrival times of the $B$ bins in increasing order (up to their first accumulation point, if they do accumulate) so that $T_1<T_2<\dots$ is the resulting sequence. The distribution of
$$\{I_m = (N_1(T_m),N_2(T_m),\dots,N_B(T_m))\}_{m\in\N}$$
is the same as that of a balls-in-bins process with feedback
function $f$ and initial conditions
$(a_1,a_2,\dots,a_B)$.\end{theorem}

One can prove this result\footnote{The exact attribution of this
result is somewhat confusing. Ref. \cite{Khanin01} cites the work
of Davis \cite{Davis90} on reinforced random walks, where it is in
turn attributed to Rubin.} as follows. First, notice that the {\em
first arrival time $T_1$} is the minimum of $X(j,a_j)$, ($1\leq
j\leq B$). By the minimum property presented above, the
probability that bin $i$ is the one at which the arrival happens
is like the first arrival probability in the corresponding
balls-in-bins process with feedback:
\begin{equation}\label{eq:min}\Pr{X(i,a_i) = \min_{1\leq j\leq
B}X(j,a_j)} = \frac{f(a_i)}{\sum_{j=1}^B f(a_j)}.\end{equation}
More generally, let $t\in\R^+$ and condition on $(N_i(t))_{i=1}^B
= (b_i)_{i=1}^B \in\N^B$, with $b_i\geq a_i$ for each $i$ (in
which case the process has not blown up). This amounts to
conditioning on
$$\forall i\in[B]\;\; \sum_{j=a_i}^{b_i-1}X(i,b_i)\leq t < \sum_{j=a_i}^{b_i}X(i,b_i).$$
From the lack of memory property of exponentials, one can deduce
that the first arrival after time $t$ at a given bin $i$ will
happen at a $\exp(f(b_j))$-distributed time, independently for
different bins. This almost takes us back to the situation of
\eqnref{min}, with $b_i$ replacing $a_i$, and we can similarly
deduce that bin $i$ gets the next ball with the desired
probability,
$$\frac{f(b_i)}{\sum_{j=1}^Bf(b_j)}.$$

\subsection{On the three regimes}\label{sec:exp_proof}

Let us now briefly point out some of the key steps in the proof of
\thmref{Khanin} via the exponential embedding, in the case $B=2$.
Reading this sketch might help the reader to become acquainted
with an important part of our methods.

We use the same notation and random variables introduced above.

Assume we start the process from state $(x,y)\in\N^2$. First, we
note that
\begin{equation}\label{eq:monopoly_cond}\Ex{\sum_{j=x}^{+\infty}X(1,j)}=\sum_{j=x}^{+\infty}
\frac{1}{f(j)},\end{equation} hence if the RHS is finite,
$\sum_{j=x}^{+\infty}X(1,j)<+\infty$ almost surely, and similarly
for $\sum_{j=y}^{+\infty}X(2,j)$. Moreover, the two random series
are independent, and neither has point-masses in their
distribution. Therefore, with probability $1$,
\begin{equation}\label{eq:dicmon}\mbox{either
}\sum_{j=x}^{+\infty}X(1,j)<\sum_{j=y}^{+\infty}X(2,j) \mbox{ or }
\sum_{j=x}^{+\infty}X(1,j)>\sum_{j=y}^{+\infty}X(2,j).\end{equation}
If the first alternative holds, there exists a finite $M>y$ such
that $$\sum_{j=x}^{+\infty}X(1,j)<\sum_{j=y}^{M-1}X(2,j).$$ Now
notice that the sequence
$\{T_{n_k}=\sum_{j=x}^{x+k}X(1,j)\}_{k\in\N}$ is an infinite
subsequence of the ball arrival times $\{T_n\}_{n\in\N}$, and at
those times
$$\forall k\in\N\;\; T_{n_k} =\sum_{j=x}^{x+k}X(1,j)<\sum_{j=y}^{M-1}X(2,j)\Rightarrow I_{n_k}(2)<M.$$
Since $\{I_m(2)\}_m$ is an increasing sequence, this means that
$I_m(2)<M$ for all $m\in\N$; that is to say, bin $1$ {\em must
achieve} monopoly. On the other hand, if the second alternative in
\eqnref{dicmon} holds, the same argument shows that bin $2$ must
achieve monopoly. Thus the condition \eqnref{monopoly_cond}
implies that with probability $1$, one of the two bins achieves
monopoly. It is not too hard to prove that if
\eqnref{monopoly_cond} does {\em not} hold, then almost surely
$\sum_{j=x}^{+\infty}X(1,j)=\sum_{j=y}^{+\infty}X(2,j)=+\infty$;
in fact, it suffices to show that, for some $\rho>0$
$$\Ex{\exp(-\rho\sum_{j=x}^{x+k}X(1,j))} = \prod_{j=x}^{x+k}\Ex{\exp(-\rho X(1,j))}\to 0\mbox{ as }k\to +\infty$$ and similarly for $\sum_{j=y}^{+\infty}X(2,j)$. In
this case one can show, by reversing the above reasoning, with
probability $1$ no bin will achieve monopoly.

Now assume that $x>y$ (for simplicity) and
\begin{equation}\label{eq:elead_cond}\sum_{j=1}^{+\infty}\frac{1}{f(j)^2}<+\infty.\end{equation}
In this case, even if $\sum_{j}f(j)^{-1}=+\infty$, the series
$$\sum_{j=x}^{+\infty} (X(1,j)-X(2,j))$$
is made of independent, centered random variables whose variances
satisfy
$$\sum_{j=x}^{+\infty}\Var{(X(1,j)-X(2,j))} = \sum_{j=1}^{+\infty}\frac{2}{f(j)^2}<+\infty.$$
Hence Kolmogorov's Three Series Theorem implies that
$\sum_{j=x}^{+\infty} (X(1,j)-X(2,j))$ converges. Following the
reasoning developed above, we deduce that almost surely
$$\mbox{either }\sum_{j=x}^{+\infty} (X(1,j)-X(2,j))-\sum_{j=y}^{x-1}X(2,j)<0\mbox{ or }\sum_{j=x}^{+\infty} (X(1,j)-X(2,j))-\sum_{j=y}^{x-1}X(2,j)>0.$$
In the first case, for all large enough $M$
$$\sum_{j=x}^{M-1}X(1,j)<\sum_{j=y}^{M-1} X(2,j)$$
and one can check that this means that for all large enough $M$,
bin $1$ reaches level $M$ {\em before} bin $2$ does (in the
embedded and continuous-time processes): that is, bin $1$ achieves
eventual leadership. Otherwise, if $\sum_{j=x}^{+\infty}
(X(1,j)-X(2,j))-\sum_{j=y}^{x-1}X(2,j)>0$, bin $2$ is the one that
achieves eventual leadership. In either case, what we have
discussed up to now proves items $1.$ and $2.$ of \thmref{Khanin}.

Finally, if
\begin{equation}\label{eq:almost_cond}\sum_{j=1}^{+\infty}\frac{1}{f(j)^2}=+\infty,\end{equation}
then for any $x$, as $k\to +\infty$,
$$\sum_{j=x}^{x+k}\frac{1}{f(j)^3}\leq \frac{1}{\min_{j\geq 1}f(j)}\, \sum_{j=x}^{x+k}\frac{1}{f(j)^2}\ll \left(\sum_{j=x}^{x+k}\frac{1}{f(j)^2}\right)^{3/2}.$$
Checking the moments of the $X(i,j)$'s and using the results in
\secref{apply_brown} shows that the sums $\sum_{j=x}^{x+k}X(1,j)$,
$\sum_{j=y}^{y+k}X(2,j)$ $(k\in\N)$ are in the domain of
attraction of Brownian Motion for any $x$ and $y$. This implies
that there is a sequence of random numbers $M_1<M_2<M_3<\dots$ and
a constant $0<\alpha<1$ such that for all $n\in\N$
$$\alpha<\Pr{\left.\sum_{j=x}^{k_n}X(1,j)<\sum_{j=y}^{k_n}X(2,j)\right| \{k_1,\dots,k_{n-1}\}\cup \{X(1,\ell),X(2,\ell):\ell\leq k_{n-1}\}}<1-\alpha.$$
This implies that both
$\sum_{j=x}^{k_n}X(1,j)<\sum_{j=y}^{k_n}X(2,j)$ and
$\sum_{j=x}^{k_n}X(1,j)>\sum_{j=y}^{k_n}X(2,j)$ must occur
infinitely often almost surely. In this case, there are infinitely
many $k$ for which bin $1$ reaches level $k$ before bin $2$ does,
and vice-versa. It follows that \eqnref{almost_cond} implies that
with probability $1$ neither bin will achieve eventual leadership,
and this proves $3.$ and the theorem.

\begin{remark}\label{rem:fictitious}Assume that bin $1$ achieves monopoly. Then all arrivals of the {\em continuous-time} process at bin $2$ after time $\sum_{j=x}^{+\infty}X(1,j)$ do not actually happen in the embedded {\em discrete-time} process $\{I_m=(I_m(1),I_m(2))\}$. We call these ``ghost events" a {\em fictitious continuation} of our process. This very useful device is akin to the continuation of a Galton-Watson process beyond its extinction time (see e.g. \cite{AlonSpencer_Method}) and is equally useful in calculations and proofs.\end{remark}

\section{Assumptions on feedback
functions}\label{sec:technical}

The purpose of this rather technical section is two-fold. First,
we spell out the technical assumptions on the feedback function
$f$ that we need in our proofs. Nothing seems to actually {\em
require} these assumptions, but they facilitate certain estimates
that we employ in the proofs.

Some readers might wish to skip the proofs in this section on a
first reading.

\subsection{Valid feedback functions}\label{sec:valid}

The feedback functions we allow in our results satisfy the
following definition.

\begin{definition}\label{def:est_smooth} A function $f:\N\to (0,+\infty)$ with $f(1)=1$\footnote{The requirement that $f(1)=1$ is just a normalization condition,
as it does not change the process.} is said to be a {\em valid
feedback function} if it can be extended to a piecewise $C^{1}$
function $g:\R^+\cup\{0\}\to(0,+\infty)$ with the following
property: if $(\ln g(\cdot))'$ is the right-derivative of $\ln g$,
and $h(x)\equiv x(\ln g(x))'$ (for $x\in\R^+\cup\{0\}$),
    \begin{enumerate}
        \item $\liminf_{x\to +\infty}h(x) = h_{min} >0$;
        \item $\lim_{x\to +\infty}x^{-{1/4}}h(x) = 0$;
        \item there exist $C>0$ and $x_0\in\R^+$ such that for all $\eps\in(0,1)$ and all $x\geq x_0$
\begin{equation}\label{eq:est_hslow}\sup_{x\leq t\leq
x^{1+\eps}}\left|\frac{h(t)}{h(x)}-1\right|\leq
C\,\eps\periodeq\end{equation} If in addition $h_{min}>1/2$, then
we say that $f$ is ELM (ELM stands for ``eventual leadership or
monopoly"). If on the other hand $h(x)\leq 1/2$ for all large
enough $x$, we say $f$ is AB (``almost-balanced").
    \end{enumerate}
With slight abuse of notation, we will always assume that $f$ is
defined over $\R^+\cup\{0\}$ and is piecewise $C^1$. We will also
call $h$ the {\em characteristic exponent} of $f$.\end{definition}
Functions with exponential growth (such as $f(x)=2^x$) or with
oscillations fail to satisfy
\defref{est_smooth}. On the other
hand, requiring that $f$ be increasing seems natural, and the
smoothness assumption still leaves us with plenty of interesting
examples of feedback functions; some examples are given in Table
\ref{tab:scaling} The ``canonical case'' where $f(x)=x^p$ ($x\geq
1$) explains the terminology for the characteristic exponent: in
that case, $h(x)\equiv p$ for all $x>1$.

\subsection{Consequences of the definition}

Let us now define the quantity
\begin{equation}S_r(n,m) \equiv \sum_{j=n}^{m-1}\frac{1}{f(n)^r}\;\;\;(r\in\N\backslash\{0\};\,n\in\N, m\in\N\cup\{+\infty\})\end{equation}
for some $f:\N\to (0,+\infty)$, and also let $S_r(n)\equiv
S_r(n,+\infty)$. If $f(x)=x^p$, then for $m-n,n\gg 1$ a simple
shows that
$$S_r(n,m)\sim\int_{n}^{m} \frac{dx}{f(x)^r} = \frac{n^{1-rp} - m^{1-rp}}{(rp-1)}.$$
The main content of the following lemmata (the first one proven in
\cite{Oliveira05}) is that a similar result holds for any valid
$f$, if $p$ is replaced by the characteristic exponent $h$. In
particular, any valid $f$ satisfies the monopoly condition in
\thmref{Khanin}. These lemmas are used in the two main proofs in
the paper.
\begin{lemma}[\cite{Oliveira05}]\label{lem:est_moment}Assume that $r$ is an integer and $f$ is a valid feedback function with characteristic exponent $h$ satisfying $h_{\min}>1/r$. Define
$$M_r(n) = \int_n^{+\infty} \frac{dx}{f(x)^r}\;\;\;\;\;(r\in\N\backslash\{0\},n\in\N)\periodeq$$
Then, as $n\to +\infty$
$$S_r(n) \sim M_r(n) \sim \frac{n}{(rh(n)-1)f(n)^r}\periodeq$$\end{lemma}

\section{Approximation by Brownian Motion}

\subsection{The Invariance Principle -- setup}\label{sec:brown}
This is the last section in which technical preliminaries are
discussed. In it, we review a form of Donsker's Invariance
Principle that shows that under suitable normalization, ``nice"
partial sums of random variables are close to Brownian Motion. All
results in this section are quite standard and can be found in
many books on Brownian Motion, e.g. \cite{BillingsleyBook}

Consider the vector space $C=C([0,1],\R)$ of all real-valued
continuous functions on the unit interval, with the sup norm
$$\|\phi(\cdot)\|_{\sup}\equiv \sup_{0\leq s\leq 1}|\phi(s)|,\;\;\;\phi(\cdot)\in C.$$
This gives $C$ a metric and a topology, and from now on we shall
think of $C$ as a measurable space with the Borel $\sigma$-field.
Brownian Motion is simply a probability measure on this measurable
space, or rather a random variable $B(\cdot)$ taking values on
$C$, whose defining properties are:

\begin{itemize}
\item $\Pr{B(0)=0} = 1$; \item for all $0\leq s_0<s_1< \dots<
s_k\leq 1$, the random variables
$$\left\{\frac{B(s_i)-B(s_{i-1})}{\sqrt{s_i-s_{i-1}}}\right\}_{i=1}^k$$ are i.i.d. standard Gaussians.
\end{itemize}

We will also use the following distributional equalities below:
\begin{eqnarray} \max_{0\leq s\leq 1}B(s),\, -\min_{0\leq t\leq 1}B(t)\eqdist |N|\mbox{ where $N$ is standard Gaussian;} \\
\label{eq:diff}\mbox{if $B'(\cdot)$ is an independent copy of
$B$}, \frac{B(\cdot)-B'(\cdot)}{\sqrt{2}} \eqdist B.\end{eqnarray}

Now consider a ``triangular sequence" $\{\xi_{n,t}\}_{t\in\N,1\leq
n\leq M_t}$ of independent, $0$-mean, square-integrable random
variables,. Letting
\begin{equation}\label{eq:scaling_sigma}\sigma^2_{k,t} \equiv
\Var{\sum_{j=1}^{k}\xi_{j,t}} =
\sum_{j=1}^{k}\Var{\xi_{j,t}}\;(k\in[M_t])\mbox{ and
}\sigma_{0,t}^2=0\commaeq\end{equation} we define a random element
$\Xi_t(\cdot)$ ($t\in\N$) of $C$ as follows.
\begin{eqnarray}\Xi_t(s) &\equiv &\frac{\sum_{j=1}^{k(s)}\xi_{j,t} + \left(\frac{s-\sigma_{k(t),m}^2}{\Var{\xi_{k(s)+1,t}}} \right)\xi_{k(s)+1,t}}{\sigma_{M_t,t}}\commaeq \\
& & \mbox{where $s\in[0,1]$ and} \\
& & k(s) \equiv \max\left\{k\in[M_t]\cup\{0\}\suchthat
\frac{\sigma_{k,t}^2}{\sigma_{M_t,t}^2}\leq
t\right\}\periodeq\end{eqnarray} Thus
$\Xi_t(\sigma_{k,t}^2/\sigma_{M_t,t}^2)$ is the sum of the $k$
first $\xi_{j,t}$'s, divided by a normalizing factor; and for
$s\in[\sigma_{k,t}^2/\sigma_{M_t,t}^2,\sigma_{k+1,t}^2/\sigma_{M_t,t}^2]$,
$\Xi_t(s)$ is defined by linear interpolation of the values of
$\Xi_t(\sigma_{k,M_t}^2/\sigma_{t,M_t}^2)$ and
$\Xi_t(\sigma_{k+1,M_t}^2/\sigma_{t,M_t}^2)$. One can check that
this is indeed a measurable element of $C$.

The Invariance Principle states that if the sequence
$\{\xi_{n,M_t}\}$ satisfies certain conditions, the distribution
of the $\Xi_t(\cdot)$'s {\em converges weakly} to a standard
Brownian motion $B(\cdot)$. What this means is that if $A\subset
C$ is measurable with boundary $\partial A$ and
$\Pr{B(\cdot)\in\partial A}=0$, then $\Pr{\Xi_t(\cdot)\in A}\to
\Pr{B(\cdot)\in A}$ as $t\to +\infty$. A sufficient condition for
this is given by

\begin{theorem}[Special case of Donsker's Invariance Principle]\label{thm:invariance}If
$$\frac{\sum_{n=1}^{M_t} \Ex{|\xi_{n,t}|^3}}{\sigma_{M_t,t}^{3/2}}\to 0 \mbox{ as }t\to +\infty,$$
then the sequence $\Xi_t(\cdot)$ converges weakly to
$B(\cdot)$.\end{theorem}

\subsection{Application to the continuous-time process}\label{sec:apply_brown}
Our typical application of this invariance principle will be to
the random variables in the exponential embedding. In the notation
of the \secref{exp_embed}, let $i,i'\in [B]$ be fixed,
$\{x_t\}_t,\{M_t\}_t\subset \N$ be sequences, and consider the
triangular array of random variables
\begin{equation}\label{eq:ourtriangle}\left\{\xi_{n,t}\equiv X(i,x_t+n-1) -
 X(i',x_t+n-1)\right\}_{m\in\N,1\leq n\leq
M_t}.\end{equation} In this case, $\sigma_{M_t,t}^2 =
2S_2(x_t,x_t+M_t)$ and the condition in \thmref{invariance} can be
seen to be equivalent to
\begin{equation}\label{eq:check}\lim_{m\to
+\infty}\frac{S_3(x_t,x_t+M_t)}{(S_2(x_t,x_t+M_t))^{3/2}}=0.\end{equation}
That is, equation \eqnref{check} is the only condition we have to
check in order to apply the Invariance Principle to the terms in
\eqnref{ourtriangle}. Notice also (by a simple limiting argument)
that if $S_2(1)<+\infty$, we can also take $M_t=+\infty$ in the
above.

\section{The scaling result for leadership}\label{sec:scaling}

\subsection{The general statement}

Recalling \secref{valid}, let $f$ be a ELM function. There exists
an $x_0$ such that $h(x) = xf'(x)/f(x)>1/2$ for all $x\geq x_0$,
which means that for such large $x$
\begin{equation}\label{eq:q0}q_0(x)\equiv \sqrt{\frac{x}{4h(x)-2}}\;\;(x\geq
x_0)\end{equation} is a well-defined, positive function. Our
generalization of \thmref{scaling_cor} shows that the quantity
$q_0(t)$ plays the same role that as the map $t\mapsto
\sqrt{t/(4p-2)}$ in the specific case $f(x)=x^p$, $p>1/2$. That
is, in order to bias a balls in bins process started with $t$
balls towards leadership by a given bin, the difference between
the initial numbers of balls should be $\theita{q_0(t)}$.

\begin{theorem}\label{thm:scaling_main}Let $f$ be a ELM function and define $q_0$ as above. Let $\lambda\in\R$ be a constant, and assume that $q=q(n)$ ($n\in\N$) is such that
\begin{itemize}
\item $t/2 \pm \lambda \,q(n)\in \N$ for all $n\in\N$; \item
$q(n)\sim q_0(n)$ for $n\gg 1$.
\end{itemize}
Now consider the $2$-bin balls-in-bins process started from
initial state
$$(x(t),y(t))=\left(\frac{t}{2}+\lambda\, q(t),\frac{t}{2}-\lambda\, q(t)\right) \in \N^2\periodeq$$
Then
\begin{eqnarray}\lim_{t\to +\infty}\Prp{(x(t),y(t))}{\elead_1} &=& \Phi(\lambda)\commaeq\\
\lim_{t\to +\infty}\Prp{(x(t),y(t))}{\lead_1} &=&
\Gamma(\lambda).\end{eqnarray}\ignore{ More generally, suppose
that we have $B\geq 2$ bins, that
$\lambda_1,\dots,\lambda_B\in\R$, are constants, and that
$q_1(\cdot),\dots,q_B(\cdot)$ are functions such that for all
$i\in\B$
\begin{itemize} \item $t/B + \lambda_i \,q_i(n)\in \N$ for all
$n\in\N$; \item $q_i(n)\sim q_0(n)$ for $n\gg 1$.\end{itemize}
Then if we start the process from state
$$(x_i(t))_{i=1}^B = \left(\frac{t}{B} + \lambda_i \,q_i(n)\right)_{i=1}^B,$$
we have
\begin{eqnarray}\lim_{t\to +\infty}\Prp{(x_i(t))_{i=1}^B}{\elead_1} & = & \Pr{\forall 2\leq i\leq B\, , N_1 - N_i < \frac{\lambda_1 - \lambda_i}{2}}\\
\lim_{t\to +\infty}\Prp{(x_i(t))_{i=1}^B}{\lead_1} & = &
\Pr{\forall 2\leq i\leq B\, , \sup_{0\leq t\leq 1}(B_1(t) -
B_i(t))< \frac{\lambda_1 - \lambda_i}{2}},\end{eqnarray} where
$\{B_i(\cdot)\}_{i=1}^B$ ($\{N_i\}_{i=1}^B$) are i.i.d. standard
Brownian Motions (resp. Gaussians).}\end{theorem}

Table \ref{tab:scaling} presents estimates of $q_0(n)$ for $n$
large, for several choices of feedback functions in the ELM
regime. In particular, the case $f(x)=x^p$ ($p>1/2$) of the
Theorem implies \thmref{scaling_cor}. The remainder of this
section contains the proof of the general result.\\
\begin{table}
\begin{center}\begin{tabular}[b]{|c|c|c|c|} \hline
  $f(x)=$ & $h(x)\sim $& $q_0(x)\sim$ & conditions \\
  \hline
  \hline
 $x^p\ln^{q}(x+e-1)$& $p$ & $\sqrt{\frac{x}{4p-2}}$ & $p>1/2$, $q\in\R$\\
\hline
$x^{q\ln^\alpha(x)}$ & $(\alpha+1)q\ln^{\alpha}(x)$ & $\sqrt{\frac{x}{4(\alpha+1)q\ln^{\alpha}(x)}}$ & $p>1/2$, $q,\alpha>0$\\
\hline $e^{x^{p}}$ & $px^p$ &
$\frac{x^{\frac{1-p}{2}}}{2\sqrt{p}}$ &
$0<p<1/4$\\
\hline
\end{tabular}\end{center}\caption{\label{tab:scaling}\small{Some
ELM feedback functions and their corresponding $h$ and $q_0$. The
last condition describes the conditions on the parameters
$p,q,\alpha$ under which each $f$ is indeed ELM.}}\end{table}

\begin{proof}[of \thmref{scaling_main}] We start
by discussing how one can write the event $\elead_1$ and lead in
terms of the exponential embedding. We only prove the result for
$\lambda>0$: the case $\lambda=0$ is a simple extension, and the
case $\lambda<0$ reduces to the one we discuss below.

Let $(x,y)\in \N^B$ be the (for the time being arbitrary) initial
conditions with $x>y$. The event $\elead_1$ holds whenever bin $1$
reaches reach level $M$ before bin $2$ does for all large enough
$M$. This requires that the time it takes for bin $1$ to reach
level $M$ in the continuous-time process is smaller than the
corresponding time for bin $2$. In the exponential embedding, this
corresponds to
$$\exists M_0\,\forall M\geq M_0\,\, \sum_{j=x}^{M}X(1,j) < \sum_{j=y}^{M} X(2,j).$$
The above event can be rewritten as
$$\exists M_0\,\forall M\geq M_0\,\, \sum_{j=x}^{M}(X(1,j) - X(2,j))< \sum_{j=y}^{x-1} X(1,j).$$
As noted in the proof sketch for \thmref{Khanin} in
\secref{exp_proof}, $\sum_{j=x}^{M}(X(1,j) - X(2,j))$ converges as
$M\to +\infty$. If follows that, except for a null event, the
above holds if and only if
$$\sum_{j=x}^{+\infty}(X(1,j) - X(2,j))< \sum_{j=y}^{x-1} X(2,j).$$
Thus we deduce that
\begin{equation}\label{eq:rewriteelead}\Prp{(x,y)}{\elead_1} = \Pr{\sum_{j=x}^{+\infty}(X(1,j) - X(2,j))< \sum_{j=y}^{x-1} X(2,j)}.\end{equation}
What about the probability of $\lead_1$? Using the above notation,
$\lead_1$ holds if for {\em all} $M\geq x$ bin $1$ reaches level
$M$ before bin $2$ does. That corresponds to
$$\forall M\geq x\,\, \sum_{j=x}^{M}X(1,j) < \sum_{j=y}^{M} X(2,j),$$
or
$$\forall M\geq x_{\ell_i}\,\, \sum_{j=x}^{M}(X(1,j) - X(2,j))< \sum_{j=y}^{x-1} X(2,j).$$
It follows that
\begin{equation}\label{eq:rewritelead}\Prp{(x,y)}{\lead_1} = \Pr{\sup_{M\geq x}\sum_{j=x}^{M}(X(1,j) - X(2,j))< \sum_{j=y}^{x-1} X(2,j)}.\end{equation}
Recall now the choice $x=x(t) = t/2 + \lambda q(t)$ and $y=y(t) =
t/2-\lambda q(t)$, where
\begin{equation}q(n) \sim \sqrt{\frac{n}{4h(n)-2}}\mbox{ for }n\gg 1\periodeq\end{equation}
As discussed in \secref{apply_brown},
$\{X(1,x(t)+n-1)-X(2,x(t)+n-1)\}_{t\in\N,n\in\N}$ is a
doubly-infinite array of centered, square-integrable random
variables with
$$\sigma_{n,t}^2\equiv \sum_{j=1}^t\Var{X(1,x(t)+n-1)-X(2,x(t)+n-1)} = 2S_2(x(t),x(t)+n) \;\; (1\leq n\leq +\infty).$$
One can construct a random continuous path $\Xi_t(\cdot)$ defined
by setting
\begin{equation}\Xi_t\left(\frac{S_2(x(t),x(t)+n)}{S_2(x(t))}\right) = \frac{\sum_{j=x}^{x(t)+n-1}X(1,j)-X(2,j)}{\sqrt{2S_2(x)}}, (1\leq n\leq +\infty)\end{equation}
and completing the remaining values by linear interpolation. As
discussed in \secref{apply_brown}, the fact that we have infinite
terms here poses no problems. Moreover, it is easy to check that
\begin{eqnarray}\sup_{0\leq s\leq 1}\Xi_t(s) &=& \frac{1}{\sqrt{2S_2(x)}}\sup_{M\geq x}\sum_{j=x}^{M}{X(1,j)-X(2,j)},\\
\Xi_t(1) &=&
\frac{1}{\sqrt{2S_2(x)}}\sum_{j=x}^{+\infty}{X(1,j)-X(2,j)}.\end{eqnarray}
Thus one can rewrite
\begin{eqnarray}\label{eq:rewriteelead2}\Prp{(x,y)}{\elead_1} &=& \Pr{\Xi_t(1) < \frac{\sum_{j=y}^{x-1} X(2,j)}{\sqrt{2S_2(x)}}},\\
\label{eq:rewritelead2}\Prp{(x,y)}{\lead_1} &=& \Pr{\sup_{0\leq
s\leq
1}\Xi_t(s)<\frac{\sum_{j=y}^{x-1}X(2,j)}{\sqrt{2S_2(x)}}}.\end{eqnarray}

Suppose we show that as $t\to +\infty$
\begin{eqnarray}\label{eq:final1}\Xi_t(\cdot) &\weakto& \mbox{a standard Brownian Motion }B(\cdot),\\
\label{eq:final2}\frac{\sum_{j=y}^{x-1}X(2,j)}{\sqrt{2S_2(x)}}&\weakto
&\lambda.\end{eqnarray} Since $\sum_{j=y}^{x-1}X(2,j)$ is
independent of $\Xi_t(\cdot)$, this means that
$$(\Xi_t(\cdot),\sum_{j=y}^{x-1}X(2,j))\weakto  (B(\cdot),\lambda),$$ which implies
\begin{eqnarray}\label{eq:limitelead2}\Prp{(x,y)}{\elead_1} &\to & \Pr{B(1)<\lambda},\\
\label{eq:limitlead2}\Prp{(x,y)}{\lead_1} &\to& \Pr{\sup_{0\leq
s\leq 1}B(s) <\lambda}.\end{eqnarray} These probabilities can be
evaluated via standard formulae for Brownian motion in
\secref{brown}, yielding the final result. We thus concentrate on
proving equations \eqnref{final1} and \eqnref{final2}.\\

{\em Proof of \eqnref{final1}.} By \secref{apply_brown}, it
suffices to show that $S_3(x(t))\ll S_2(x(t))^{3/2}$. But this
follows directly from the formulae in \lemref{est_moment}, the
fact that $x=x(t)\to +\infty$, and the assumption that $h(x)\ll
\sqrt{x}$:
$$S_3(x) \sim \frac{x}{(3h(x)-1)f(x)}\ll \left(\frac{x}{(2h(x)-1)f^2(x)}\right)^{3/2} = S_2(x)^{3/2}.$$\\

{\em Proof of \eqnref{final2}.} Let us first establish a few facts
about $f$, $x$, $y$ and $q$.
\begin{enumerate}
\item {\em For $t\gg 1$, $q(t) \sim
f(t/2)\sqrt{S_2(t/2)/2}=\bigoh{t}$}. Indeed, for $t$ large,
\lemref{est_moment} implies that $S_2(t/2)\sim
t/(4h(t/2)-2)f(t/2)^2$, and, because of \eqnref{est_hslow} in
\defref{est_smooth}, $h(t/2)\sim h(t)$ . Moreover, since $\liminf_{n\to
+\infty}h(n)>1/2$, $q(t) = \sqrt{t/\ohmega{1}} =
\bigoh{\sqrt{t}}$. \item {\em For $t\gg 1$, $f(x)=f(t/2 + \lambda
q(t))\sim f(t/2)$}. In this case
\begin{eqnarray}\left|\ln \frac{f(x)}{f\left(\frac{t}{2}\right)}\right| &=&
\left|\int_{\frac{t}{2}}^{\frac{t}{2} + \lambda q(t)} (\ln
f(u))'\,du\right| \\ \nonumber &=&
\left|\int_{\frac{t}{2}}^{\frac{t}{2} + \lambda q(t)}
\frac{h(u)}{u}\,du \right|\\ \nonumber & \leq &
\sup\limits_{\frac{t}{2}\leq u\leq \frac{t}{2} + \lambda q(t)}
|h(u)|\, \ln\left(1 + \frac{2
q(t)}{t}\right)\periodeq\end{eqnarray} Now notice that equation
\eqnref{est_hslow} implies that the $\sup$ of $h(u)$ above is
$\sim h(t/2)\sim h(t)$. Moreover, $q(t)/t = \sqrt{1/t(2h(t)-1)}\ll
1$, since $\liminf_{n\to +\infty}h(n)>1/2$. Therefore,
\begin{equation}\left|\ln
\frac{f(x)}{f\left(\frac{t}{2}\right)}\right| =
\bigoh{\frac{h(t)q(t)}{t}} = \bigoh{\frac{1}{\sqrt{t}}} =
\liloh{1}\periodeq\end{equation} \item {\em For $t\gg 1$,
$f(y)\sim f(t/2)$}. The proof is almost identical to the one
above. \item {\em For $t\gg 1$ and $r\geq 2$, $S_r(x),S_r(y)\sim
S_r(t/2)$}. Indeed, because $q(t)=\bigoh{\sqrt{t}}$, $x,y\to
+\infty$,, and the formulae in \lemref{est_moment} apply. Thus
$S_2(x)\sim x/[(rh(x)-1)f(x)^2]$, and by $1.$ and $3.$, $x\sim
t/2$, $h(x) \sim h(t/2) $ and $f(x)\sim f(t/2)$, which implies
$S_2(x) \sim (t/2)/[(rh(t/2)-1)f(t/2)^2]\sim S_r(t/2)$. The same
argument proves the desired result for $S_r(y)$.\end{enumerate}

We now apply the estimates to the problem at hand. For $t$ large
enough (so that $x$ and $y$ are also large), we can ensure that
$f$ is increasing on $[y,x]$, so that
$$\frac{x-y}{f(x)} \leq \Ex{\sum_{\ell=y}^{x-1} X(2,\ell)} = \sum_{j=y}^{x-1}\frac{1}{f(j)} \leq
\frac{x-y}{f(y)}\periodeq$$ By items $2.$ and $3.$ above and the
definition of $x,y$, this implies that
\begin{equation}\Ex{\sum_{\ell=y}^{x-1} X(2,\ell)} \sim \frac{2\lambda
q(t)}{f(t/2)}\;\;(t\gg 1)\periodeq\end{equation} Similarly, one
can show that
\begin{equation}\Var{\sum_{\ell=y}^{x-1} X(2,\ell)} \sim \frac{2\lambda
q(t)}{f(t/2)^2}\;\;(t\gg 1)\periodeq\end{equation} Therefore,
using $4.$,
\begin{equation}\Ex{\sum_{\ell=y}^{x-1} X(2,\ell)}^2 \geq
(1-\liloh{1})2\lambda\,q(t)\,\Var{\sum_{\ell=y}^{x-1}
X(2,\ell)}\;\;(t\gg 1)\periodeq\end{equation} By Chebyshev's
Inequality, it follows that
$$\frac{{\sum_{\ell=y}^{x-1} X(2,\ell)}}{\left(\frac{2\lambda
q(t)}{f(t/2)}\right)}\weakto 1\mbox{ as $t\to +\infty$.}$$
Finally, since
$$\frac{2q(t)}{f(t/2)}\sim \sqrt{\frac{2t}{(2h(t)-1)f(t/2)^2}}\sim \sqrt{2S_2(x)},$$
we have
$$\frac{{\sum_{\ell=y}^{x-1} X(2,\ell)}}{\lambda
\,\sqrt{2S_2(x)}}\weakto 1\mbox{ as $t\to +\infty$,}$$ which is
the desired result.
\end{proof}

\section{The almost balanced regime}\label{sec:almost}

This section proves our result on the overtaking time, a
generalization of \thmref{almost_cor}. To recapitulate:
\thmref{Khanin} tells us that, when the feedback function $f$
satisfies
\begin{equation}\label{eq:almost_ab}\sum_{j=1}^{+\infty}\frac{1}{f(j)^2}=+\infty\commaeq\end{equation}
each of the two bins will be the one with more balls infinitely
many times. Our main interest in this chapter will be in
determining {\em how long} it takes for bin $1$ to have more balls
than bin $2$, given that the latter bin has more balls at the
start. More specifically, assume the process
$(I_1(\cdot),I_2(\cdot))$ is started from state
\begin{equation}(\lceil \alpha t\rceil,t - \lceil \alpha
t\rceil)\commaeq\;\; t\gg 1\mbox{ and } \alpha\in(0,1/2)\mbox{
fixed}\periodeq\end{equation} As in the introduction, let $V$ be
the {\em overtaking time} of the process: that is the first time
when bin $1$ has more balls than bin $2$.
\begin{equation}\label{eq:almost_overtaking}V \equiv \min\{v\in\N\suchthat
I_1(v)>I_2(v)\}\periodeq\end{equation} Under condition
\eqnref{almost_ab}, this $\min$ exists and is finite with
probability $1$. We will be interested in describing the
asymptotic distribution of $V$.

To express our main result, let us introduce two mappings.

\begin{eqnarray}\funcdef{F_{t,\alpha}}{((1-\alpha)t,+\infty)}{\R^+}{u}{\frac{\int_{\alpha t}^{(1-\alpha)
t}\frac{dx}{f(x)}}{\sqrt{\int_{(1-\alpha)
t}^u\frac{dx}{f(x)^2}}}}\commaeq\\
G_{t,\alpha}\equiv \mbox{the inverse of
}F_{t,\alpha}\periodeq\end{eqnarray}

Notice that $\lim_{u\searrow (1-\alpha)t}F_{t,\alpha}(u)=+\infty$
and
$$\lim_{u\to +\infty}F_{t,\alpha}(u) =\frac{\int_{\alpha t}^{(1-\alpha)
t}\frac{dx}{f(x)}}{\sqrt{\int_{(1-\alpha)
t}^{+\infty}\frac{dx}{f(x)^2}}}=0$$ because (as a consequence of
$\sum_{j}^{+\infty}f(j)^{-2}=+\infty$) the denominator in the RHS
is infinite. Thus $F_{t,\alpha}$ is a monotone-decreasing function
whose range is $\R^+$, and $G_{t,\alpha}$ is not only
well-defined, but monotone-decreasing as well.

\begin{theorem}\label{thm:almost_main}Assume that $f$ is a AB function
(cf. \defref{est_smooth}), and define $F_{t,\alpha}$,
$G_{t,\alpha}$ as above. Let $V_{t,\alpha}$ be the random variable
$V$ defined above, conditioned on the initial state $(\lceil
\alpha t\rceil,t-\lceil \alpha t\rceil)$ of the balls-in-bins
process. Then, as $t\to +\infty$,
\begin{equation}\forall \lambda\in\R^+\commaeq\;\lim_{t\to +\infty}\Pr{V_{t,\alpha}\geq
2G_{t,\alpha}(\lambda) - (t+1)}=
\Gamma(\lambda)\commaeq\end{equation} (where $\Gamma(\cdot)$ is
defined as in \thmref{scaling_cor}), or equivalently,
\begin{equation}F_{t,\alpha}\left(\frac{V_{t,\alpha}+ (t+1)}{2}\right)\weakto
|N|\commaeq \end{equation} for a standard Gaussian random variable
$N$. The latter expression makes sense because $V_{t,\alpha}\geq
(1-\alpha)t$, so
$$\frac{V_{t,\alpha}+ (t+1)}{2}\geq \left(1-\frac{\alpha}{2}\right)t > (1-\alpha)t\periodeq$$\end{theorem}

This result is quite general, but applying it to a specific
situation requires a calculation. We do this for the case
$f(x)=x^p$ below, and then prove \thmref{almost_main} below.

\subsection{Proof of the special case}

\begin{claim}\thmref{almost_main} implies \thmref{almost_cor}.\end{claim}

\begin{proof} One way of interpreting \thmref{almost_main} in the $f(x)=x^p$ case is by
saying that
$$V_{t,\alpha} \equiv 2G_{t,\alpha}(U_{t,\alpha,p})-(t+1)\mbox{ with probability $\to 1$}\commaeq$$
where $U_{t,\alpha,p}\weakto |N|$ as $t\to +\infty$. Thus the
Corollary follows from providing a formula for $G_{t,\alpha}$.  We
will first assume that $0<p<1/2$, in which case
\begin{eqnarray}\label{eq:almost_eta} e(t,\alpha) &\equiv& \int_{\alpha t}^{(1-\alpha)t}\frac{dx}{f(x)}\\
&=&[(1-\alpha)^{1-p} -
\alpha^{1-p}]\frac{t^{1-p}}{1-p}\commaeq\end{eqnarray} and for all
$u> (1-\alpha)t$
\begin{eqnarray}\label{eq:almost_feta} f(u,t,\alpha) &\equiv& \int_{(1-\alpha) t}^{u}\frac{dx}{f(x)^2}\\
&=&\frac{u^{1-2p}- [(1-\alpha)
t]^{1-2p}}{1-2p}\periodeq\end{eqnarray} Then, for $u\geq
(1-\alpha)t$, \begin{eqnarray}F_{t,\alpha}(u) &=&
\frac{e(t,\alpha)}{\sqrt{f(u,t,\alpha)}} \\
&=&\frac{[(1-\alpha)^{1-p} -\alpha^{1-p}]\frac{t^{1-p}}{1-p}
}{\sqrt{\frac{u^{1-2p}- [(1-\alpha)
t]^{1-2p}}{1-2p}}}\end{eqnarray} To compute
$G_{t,\alpha}(\lambda)$ for some $\lambda\in\R^+$, we must solve
the equation $F_{t,\alpha}(G_{t,\alpha}(\lambda))=\lambda$, which
corresponds to
\begin{equation}G_{t,\alpha}(\lambda)^{1-2p} = [(1-\alpha) t]^{1-2p} + \frac{(1-2p)}{\lambda^2}\left\{[(1-\alpha)^{1-p} -\alpha^{1-p}]\frac{t^{1-p}}{1-p}\right\}^2\periodeq\end{equation}

Therefore,
\begin{equation}G_{t,\alpha}(\lambda) =  \left\{\frac{(1-2p)}{(1-p)^2}[(1-\alpha)^{1-p} -\alpha^{1-p}]^2 + \frac{(1-\alpha)^{1-2p}\lambda^2}{t}\right\}^{\frac{1}{1-2p}}\frac{t^{1+\frac{1}{1-2p}}}{\lambda^{1+\frac{2p}{1-2p}}}\commaeq\end{equation}
and the result for $0<p<1/2$ follows.

For $p=1/2$, the above formula for $e(t,\alpha)$ still applies,
but
\begin{eqnarray}\label{eq:almost_feta12} f(u,t,\alpha) &\equiv& \int_{(1-\alpha) t}^{u}\frac{dx}{f(x)^2}\\
&=&\ln\frac{u}{(1-\alpha)t}\commaeq\end{eqnarray} and thus
\begin{eqnarray}G_{t,\alpha}(\lambda) &=& (1-\alpha)\, t\,
\exp\left\{\frac{e(t,\alpha)^2}{\lambda^2}\right\} \\ &=&
(1-\alpha)\, t\,
\exp\left\{4[1-2\sqrt{\alpha(1-\alpha)}]\,\frac{t}{\lambda^2}\right\}
\periodeq\end{eqnarray} This finishes the proof. \end{proof}

\subsection{Proof of \thmref{almost_main}}\label{sec:almost_proof}
\begin{proof} Define the {\em overtaking number} $N_{t,\alpha}$ to be
$I_1(V_{t,\alpha})$, i.e. the number $I_1(v)$ of balls in bin $1$
at the first time $v$ when $I_1(v)>I_2(v)$, under initial
conditions $(\lceil \alpha t\rceil,t-\lceil \alpha t\rceil)$. It
follows from this definition that at time $v'=V_{t,\alpha}-1$,
$I_1(v')=I_2(v') = N_{t,\alpha}-1$. Since $I_1(0)+I_2(0)=t$, this
means that
\begin{eqnarray}\label{eq:almost_overtake}\nonumber V_{t,\alpha}-1 + t &=& I_1(v')+I_2(v') \\ \nonumber &=&\mbox{total $\#$ of balls at time $v'$}\\ \nonumber
&=& 2N_{t,\alpha} - 2 \\ \label{eq:almost_overtakelast}\Rightarrow
V_{t,\alpha} &=& 2N_{t,\alpha} - (t+1)\periodeq\end{eqnarray} Thus
results about the distribution of $N_{t,\alpha}$ translate
immediately into results about $V_{t,\alpha}$. Since
$N_{t,\alpha}$ is easier to analyze via our techniques, we shall
spend most of our time considering this quantity, returning to the
more significant $V_{t,\alpha}$ at the end of the proof.

We begin by showing that, in terms of the exponential embedding
random variables,
\begin{equation}\label{eq:almost_expwrite}\forall M\in\N\;\; \{N_{t,\alpha}\geq M\} =
\left\{\sup_{y\leq m\leq M-1}\sum_{j=y}^{m-1}(X(2,j)-X(1,j))\leq
\sum_{\ell=x}^{y-1}X(1,\ell)\right\}\commaeq\end{equation} where
$x\equiv\lceil\alpha t\rceil$ (respectively $y\equiv
t-\lceil\alpha t\rceil$) is the initial number of balls in bin $1$
(resp. $2$). Indeed, $N_{t,\alpha}\geq M$ occurs if and only if
for all $y\leq m\leq M-1$, the time it takes for bin $2$ to
receive its $m$th ball (which is $\sum_{j=y}^{m-1}X(2,j)$) is
smaller than or equal to the time it takes for bin $1$ to receive
its $m$th balls (which is $\sum_{\ell=x}^{m-1}X(1,\ell)$).
Symbolically,
\begin{equation*}\label{eq:almost_expwrite2}\forall M\in\N\;\; \{N_{t,\alpha}\geq M\} =
\left\{\forall y\leq m\leq M-1\commaeq\;
\sum_{j=y}^{m-1}X(2,j)\leq \sum_{\ell=x}^{m-1}
X(1,\ell)\right\}\commaeq\end{equation*} from which
\eqnref{almost_expwrite} follows.

We now wish to estimate the probability of the event at the RHS of
\eqnref{almost_expwrite}. We begin by looking at
\begin{equation}\label{eq:partialsumsalmost}\sup_{y\leq m\leq M-1}\sum_{j=y}^{m-1}(X(2,j)-X(1,j))\periodeq\end{equation}
In particular, let us assume that
\begin{equation}\label{eq:almost_lambda}M \equiv \lfloor G_{t,\alpha}(\lambda)\rfloor \end{equation}
for some $\lambda\in\R^+$, so that
\begin{equation}\{N_{t,\alpha}\geq M\} = \{N_{t,\alpha}\geq G_{t,\alpha}(\lambda)\}\periodeq\end{equation}
We wish to apply the Invariance Principle to the random variables
$$\{\xi_{n,t}\equiv X(2,y(t)+n-1)-X(1,y(t)+n-1)\}_{t\in\N,1\leq n\leq M-y(t)-1}.$$
Note that
$$\sigma^2_{n,t}\equiv \sum_{j=1}^n \Var{X(2,y(t)+n-1)-X(1,y(t)+n-1)} = 2S_2(y(t),y(t)+n).$$
Following \secref{apply_brown}, we build the random path
$\Xi_t(\cdot)$ that linearly interpolates the values
$$\Xi_t\left(\frac{S_2(y(t),y(t)+n)}{S_2(y(t),M-1)}\right) \equiv \frac{\sum_{j=1}^n {X(2,y(t)+j-1)-X(1,y(t)+j-1)}}{\sqrt{2S_2(y,M-1)}}, \, 1\leq n\leq M-y(t).$$
As in the previous proof, we have
$$\sup_{0\leq s\leq 1}\Xi_t(s) = \frac{\sup_{y\leq m\leq M-1}\sum_{j=y}^{m-1} {X(2,y(t)+j-1)-X(1,y(t)+j-1)}}{\sqrt{2S_2(y,M-1)}}.$$
Hence
\begin{equation*}\Pr{N_{t,\alpha}\geq M} = \Pr{\sup_{0\leq s\leq 1}\Xi_t(s) \leq \frac{\sum_{j=x}^{y-1}X(1,j)}{\sqrt{2S_2(y,M-1)}}}.\end{equation*}
We will eventually prove that as $t\to +\infty$,
\begin{eqnarray}\label{eq:almost_Cheby}\sum_{j=x}^{y-1}X(1,j)&\weakto& \lambda,\\
\label{eq:almost_Brown}\Xi_t(\cdot)&\weakto& \mbox{a standard
Brownian Motion }B(\cdot).\end{eqnarray} It follows from this and
the independence of $\Xi_t(\cdot)$ $\sum_{j=x}^{y-1}X(1,j)$ that
\begin{equation}\lim_{t\to +\infty}\Pr{N_{t,\alpha}\geq M} =\Pr{\sup_{0\leq s\leq 1}B(s) \leq \lambda} = \Gamma(\lambda),\end{equation}
which is the desired result. Thus we concentrate on proving
\eqnref{almost_Brown} and \eqnref{almost_Cheby}.\\

{\em Proof of \eqnref{almost_Cheby}. } The expectation of
$\sum_{\ell =x}^{y-1} X(1,\ell)$ can be estimated as follows.
\begin{eqnarray} \Ex{\sum_{\ell =x}^{y-1} X(1,\ell)} &= &S_1(\lceil \alpha t\rceil, t-\lceil \alpha t\rceil)\\
\label{eq:integral}&\sim &\int_{\alpha
t}^{(1-\alpha)t}\frac{ds}{f(s)}\\ \label{eq:ohmega}&=&
\ohmega{\frac{t}{f(\alpha t)}}\gg 1.\end{eqnarray}Indeed,
\eqnref{integral} follows from the fact that $\sum_{j\geq
1}f(j)^{-1}=+\infty$, and \eqnref{ohmega} follows from the
assumption $h(s)\leq 1/2$ for $s$ large, which means that
$f((1-\alpha)t)=\bigoh{f(\alpha t)}=\bigoh{\sqrt{t}}$. The
variance of the is
\begin{equation*}\Var{\sum_{\ell =x}^{y-1} X(1,\ell)} = S_2(\lceil\alpha t\rceil,t-\lceil\alpha
t\rceil) =\bigoh{\frac{t}{f(\alpha t)^2}} \ll \Ex{\sum_{\ell
=x}^{y-1} X(1,\ell)}^2 \;\;(t\gg 1)\end{equation*} because for $s$
large, $h(s)\geq 0$ and thus $f$ is decreasing. By Chebyshev's
Inequality, it follows that
\begin{equation}\label{eq:convCheby}\frac{\sum_{\ell =x}^{y-1} X(1,\ell)}{\int_{\alpha
t}^{(1-\alpha)t}\frac{ds}{f(s)}}\weakto 1.\end{equation}

On the other hand, notice that if $M'\equiv G_{t,\alpha}(\lambda)$
satisfies
\begin{equation}\label{eq:definitionFok}F_{t,\alpha}(M') =
\frac{\int_{\alpha
t}^{(1-\alpha)t}\frac{ds}{f(s)}}{\sqrt{2\int_{(1-\alpha)
t}^{M'}\frac{ds}{f(s)^2}}} = \lambda\periodeq\end{equation} We
wish to show that
\begin{equation}\label{eq:simisgood}F_{t,\alpha}(M)\sim
\lambda;\end{equation} this will follow from
$$\int_{(1-\alpha)t}^{M'}\frac{ds}{f(s)^2} \sim \int_{(1-\alpha)
t}^{M}\frac{ds}{f(s)^2}.$$ Since $M=\lfloor M'\rfloor$, $M\leq
M'\leq M+1$ and
$$\left|\int_{(1-\alpha) t}^{M'}\frac{ds}{f(s)^2}-\int_{(1-\alpha)
t}^{M}\frac{ds}{f(s)^2}\right| \leq \frac{1}{f(M')^2}.$$ This
follows from the fact that $h(s)>0$ (and thus $f$ is increasing)
for $s$ large enough. On the other hand, we must have
$M'-(1-\alpha)t\gg 1$, as for any constant $C$
$$\int_{(1-\alpha)t}^{(1-\alpha)t+C}\frac{ds}{f(s)^2}\leq \frac{C}{f((1-\alpha)t)}\ll 1\ll \left(\int_{\alpha
t}^{(1-\alpha)t}\frac{ds}{f(s)}\right)^2,$$ contradicting
\eqnref{definitionFok}. Thus
$$\int_{(1-\alpha)t}^{M'}\frac{ds}{f(s)^2}\geq \frac{M'-(1-\alpha)t}{f(M')^2}\gg \left|\int_{(1-\alpha) t}^{M'}\frac{ds}{f(s)^2}-\int_{(1-\alpha)
t}^{M}\frac{ds}{f(s)^2}\right|,$$ which proves \eqnref{simisgood}.
In particular, it follows from \eqnref{convCheby} that
\begin{equation}\label{eq:almost_sec}\frac{\sum_{j=x}^{y-1}X(1,j)}{\int_{(1-\alpha)
t}^{M}\frac{ds}{f(s)^2}}\weakto \lambda.\end{equation} Finally,
since by \eqnref{ohmega} the numerator of
$$F_{t,\alpha}(M) = \frac{\int_{\alpha t}^{(1\alpha)t} \frac{ds}{f(s)^2}}{\sqrt{2\int_{(1-\alpha)t}^M \frac{ds}{f(s)^2}}}$$
diverges, we deduce that $\int_{(1-\alpha)t}^M
\frac{ds}{f(s)^2}\gg 1$, which implies that
$$\int_{(1-\alpha)t}^M
\frac{ds}{f(s)^2}\sim S_2(y,M-1).$$ Plugging this into
\eqnref{almost_sec} yields \eqnref{almost_Cheby}.\\

{\em Proof of \eqnref{almost_Brown}. } We have already shown above
that $S_2(y,M-1)\sim \int_{(1-\alpha)t}^M \frac{ds}{f(s)^2}\gg 1$.
Since $f$ is bounded below,
$$S_3(y,M-1) = \bigoh{S_2(y,M-1)}\ll S_2(y,M-1)^{3/2}.$$
Thus the condition for the Invariance Principlie in
\secref{apply_brown} is satisfied, and this finishes the
proof.\end{proof}

\section{Extensions and open problems}\label{sec:last}

Our \thmref{scaling_main} admits an extension to the case of a
general number $B\geq 2$ of bins.

\begin{theorem}\label{thm:scaling_main2}Let $f$ be a ELM function, $B\geq 2$ and $\lambda_i\in\R$ ($i\in B$) be constants, and assume that $q=q(n)$ ($n\in\N$) is such that
\begin{itemize}
\item $t/B +\lambda_i \,q(n)\in \N$ for all $n\in\N$ and
$i\in[B]$; \item $q(n)\sim q_0(n)$ for $n\gg 1$ ($q_0$ is defined
in \eqnref{q0}).
\end{itemize}
Now consider the $B$-bin balls-in-bins process started from
initial state
$$(x_i(t))_{i=1}^B =\left(\frac{t}{B}+\lambda_i\, q(t)\right)_{i=1}^B \in \N^B\periodeq$$
Then
\begin{eqnarray*}\lim_{t\to +\infty}\Prp{(x_i(t))_{i=1}^B}{\elead_1} & = & \Pr{\forall 2\leq i\leq B\, , \frac{N_1-N_i}{\sqrt{2}} < (\lambda_1 - \lambda_i)}\\
\lim_{t\to +\infty}\Prp{(x_i(t))_{i=1}^B}{\lead_1} & = &
\Pr{\forall 2\leq i\leq B\, , \sup_{0\leq t\leq 1}\frac{(B_1(t) -
B_i(t))}{\sqrt{2}}< (\lambda_1 - \lambda_i)},\end{eqnarray*} where
$\{B_i(\cdot)\}_{i=1}^B$ ($\{N_i\}_{i=1}^B$) are i.i.d. standard
Brownian Motions (resp. Gaussians).\end{theorem}

The proof of this result, which we omit, follows essentially the
same lines as that of \thmref{scaling_main}. The only major
differences is that the Invariance Principle is applied to the
sequences
$$\{\xi^{(i)}_{n,t}=X(i,x_i(t)+n-1) - f(x_i(t)+n-1)^{-1}\}_{t\in\N, n\in\N}\,\,(i\in[B])$$
and some extra care must be taken in looking at the maxima of the
corresponding random continuous functions.

It is not entirely obvious how \thmref{almost_main} on the almost
balanced regime should or could be generalized to more than $2$
bins. All we can say in general is that the addition of bins to
the process will always make any overtaking take longer (to see
this, notice that more bins just add more arrivals in the
continuous-time exponential embedding between the start and
overtaking times). In fact, there is a more basic question about
such processes that we cannot answer; it was posed as a conjecture
by Joel Spencer.

\begin{conjecture}[Joel Spencer]Consider a balls-in-bins process $(I_m(i))_{i\in[B],m\geq 0}$ with feedback function $f$ that is in the almost-balanced regime (cf. \thmref{Khanin}). Then for all permutations $\Pi$ of $[B]$ and all initial conditions,
with probability $1$ there are infinitely many $m\geq 0$ with
$I_m(\Pi(1))<I_m(\Pi(2))<\dots<I_m(\Pi(B))$. That is, all
permutations possible of the bins occur infinitely often almost
surely.\end{conjecture}

Perhaps the techniques presented here could be used to settle this
problem.

\bibliography{binb_brown}

\begin{thebibliography}{10}

\bibitem{AlbertSurvey}
R\'eka Albert and Albert-L\'aszl\'o Barab\'asi.
\newblock Statistical mechanics of complex networks.
\newblock {\em Reviews of Modern Physics}, 74:47--97, 2002.
\newblock Arxiv: \texttt{cond-mat/0106096}.

\bibitem{AlonSpencer_Method}
Noga Alon and Joel Spencer.
\newblock {\em The Probabilistic Method}.
\newblock Wiley-Interscience Series in Discrete Mathematics. John Wiley and
  Sons, New York, second edition, 2000.

\bibitem{BillingsleyBook}
Patrick Billingsley.
\newblock {\em Convergence of Probability Measures}.
\newblock Wiley Series in Probability and Statistics. John Wiley and Sons, New
  York, second edition, 1999.

\bibitem{Davis90}
Burgess Davis.
\newblock Reinforced random walk.
\newblock {\em Probability Theory and Related Fields}, 84(2):203--229, 1990.

\bibitem{DrineaEM}
Eleni Drinea, Mihaela Enachescu, and Michael Mitzenmacher.
\newblock Variations on {Random Graph} models of the {Web}.
\newblock {Harvard Technical Report TR}-06-01, 2001.

\bibitem{Drinea02}
Eleni Drinea, Alan Frieze, and Michael Mitzenmacher.
\newblock Balls in bins processes with feedback.
\newblock In {\em Proceedings of the 11th Annual {ACM-SIAM} Symposium on
  Discrete Algorithms}, pages 308--315. Society for Industrial and Applied
  Mathematics, Philadelphia, PA, USA, 2002.

\bibitem{Khanin01}
Kostya Khanin and Raya Khanin.
\newblock A probabilistic model for the establishment of neuron polarity.
\newblock {\em Journal of Mathematical Biology}, 42(1):26--40, 2001.

\bibitem{Krapivsky02}
P.L. Krapivsky and Sidney~L. Redner.
\newblock Organization of growing random networks.
\newblock {\em Physics Reviews E}, 63:066123, 2001.
\newblock Available at Arxiv: \texttt{cond-mat/0011094}.

\bibitem{MitzenmacherOS04}
Michael Mitzenmacher, Roberto Oliveira, and Joel Spencer.
\newblock A scaling result for explosive processes.
\newblock {\em Electronic Journal of Combinatorics}, 11(1):R31, 2004.

\bibitem{tese}
Roberto Oliveira.
\newblock {\em Preferential attachment}.
\newblock PhD thesis, Department of Mathematics, Courant Institute of
  Mathematical Sciences, New York University, 2004.

\bibitem{Oliveira05}
Roberto Oliveira.
\newblock {The Onset of Dominance in Balls-in-Bins Processes with Feedback}.
\newblock Preprint. Arxiv: \texttt{math.PR/0510415}, 2005.

\bibitem{avoid}
Roberto Oliveira and Joel Spencer.
\newblock Avoiding an iminent defeat in a balls-in-bins process with feedback.
\newblock Manuscript, 2005.

\bibitem{Spencer??}
Joel Spencer and Nicholas Wormald.
\newblock Explosive processes.
\newblock Manuscript.

\end{thebibliography}
\bibliographystyle{plain}

\end{document}